
\input amstex 
\documentstyle{amsppt}
\input bull-ppt
\topmatter
\cvol{27}
\cvolyear{1992}
\cmonth{July}
\cyear{1992}
\cvolno{1}
\cpgs{143-147}
\title The Green function of Teichm\"uller\\
spaces with applications\endtitle
\shorttitle{Green Function of Teichm\"uller Spaces}
\author Samuel L. Krushkal\endauthor
\shortauthor{S. L. Krushkal}
\address Research Institute for Mathematical
Sciences and Department of Mathematics and Computer 
Science, Bar-Ilan 
University, 52900 Ramat-Gan, Israel\endaddress
\address Institute of Mathematics, Siberian Branch of the 
Russian Academy of 
Sciences, 630090, Novosibirsk, Russia\endaddress
\date November 1, 1991\enddate
\abstract We describe briefly a new approach to some 
problems related to 
Teichm\"uller spaces, invariant metrics, and extremal 
quasiconformal maps. 
This approach is based on the properties of 
plurisubharmonic functions,
especially of the plurisubharmonic Green function.

The main theorem gives an explicit representation of the 
Green function for
Teichm\"uller spaces by the Kobayashi-Teichm\"uller metric 
of these
spaces.  This leads to various applications.  
In particular, this gives a new characterization of 
extremal quasiconformal
maps.
\endabstract
\endtopmatter

\document
\par I want to describe here briefly a somewhat new 
approach to certain subjects
related to the Teichm\"uller theory, hyperbolic geometry, 
and quasiconformal
maps.
This approach is based on the properties of 
plurisubharmonic functions,
especially of the Green function.
\heading1. The Green function\endheading
\par
Recall that the Green function $g_D(x,y)$ of a domain $D$ 
in a complex Banach
space $\Cal X$ is defined as
$$g_D(x,y)=\sup u_y(x) \quad (x,y\in D),\tag1$$
where supremum is taken over all plurisubharmonic 
functions $u_y\:
D\rightarrow [-\infty,0)$ having  the
representation $u_y(x)=\log \|x-y\|+O(1)$
in a neighborhood of the point $y$; here $\|\boldcdot\|$ 
is the norm on
$\Cal X$, and $O(1)$ is bounded above.  
This definition is easily extended to complex Banach 
manifolds.

The important properties of the Green function are: 
\roster
\item"(i)" $g_D(x,y)$ does not increase under holomorphic 
maps from $D;$
\item"(ii)" $g_D(x,y)$ is plurisubharmonic in $x$ on $D$ 
for every fixed $y\in
D;$
\item"(iii)"for the domains $D\subset \bold C^n,$ the 
function $g_D(x,y)$
(provided it belongs to $L^{\infty,\operatorname{loc}}
(D\backslash\{y\})$ is maximal and is
a generalized solution of the homogeneous Monge-Amp\'ere 
equation
$(dd^cu)^n=0$ in $D\backslash\{y\}$.
\endroster
(For details and other properties we refer
to \cite{2, 7, 11, 14}.)
\heading2. Main theorem\endheading

Now let us consider the Teichm\"uller space $T(\Gamma)$ of 
an arbitrary
Fuchsian group $\Gamma,$ for simplicity, without torsion.
It is well known that $T(\Gamma)$ admits a structure of a 
complex (Banach)
manifold; moreover, it can be holomorphically embedded as 
a bounded domain
into some complex Banach space (see, e.g., \cite{5, 12}).

Our main result states:

\proclaim{Theorem 1}
The Green function $g_{T(\Gamma)}(x,y)$ of the space 
$T(\Gamma)$ is given by
$$g_{T(\Gamma)}(x,y)=\log\frac{e^{2d_{T(\Gamma)}(x,y)}-1}
{e^{2d_{T(\Gamma)}(x,y)}+1}\tag2$$
where $d_{T(\Gamma)}$ is the Kobyashi-Teichm\"uller metric 
on $T(\Gamma).$
\endproclaim
\heading3. Some corollaries\endheading

Taking into account the known properties of the 
Kobayashi-Teichm\"uller
metric, one derives  from (2) the following remarkable 
consequences.

\proclaim{Corollary 1}
The function $g_{T(\Gamma)}$ is symmetric with respect to 
its arguments
$$g_{T(\Gamma)}(x,y)=g_{T(\Gamma)}(y,x).$$
\endproclaim
\proclaim{Corollary 2}
For the finite-dimensional Teichm\"uller space $T(p,n)$ 
\RM(which corresponds to
the Riemann surfaces of finite conformal type $(p,n)$\RM), 
the Green function
$g_{T(\Gamma)}(x,y)$ is $C^1$-differentiable on 
$T(p,n)\times
T(p,n)\backslash\{diagonal\}.$
\endproclaim

The latter result is of special interest because the  
known general results
about the Green function $g_D$ state at most the 
continuity of $g_D$ on
$D\times D,$ even in the case of domains $D$ in $\bold 
C^n$ with
$C^\infty$-boundaries.

Corollary 2 follows from equality (2)  and the
$C^1$-differentiability of the Teichm\"uller metric for 
the spaces $T(p,n)$
established in \cite{3,\ 5}.  
On the other hand, it is well known (see \cite{5,\ 16}) that
$d_{T(\Gamma)}$ is not $C^2$-differentiable (along some 
directions in
$T(p,n))$; hence, $g_{T(\Gamma)}$ also is not 
$C^2$-differentiable in $x$ on
$T(p,n)\backslash \{y\}.$

Now let us consider an extremal Beltrami differential 
$\mu(z)d\overline z/dz$
on the Riemann surface $X$ represented by the point $x$ in 
$T(\Gamma)$ that
realizes the distance $d_{T(\Gamma)}(x,y)$; that is, 
having denoted
$k(x,y)=\|\mu\|_{L^\infty}$ we get
$$d_{T(\Gamma)}(x,y)=\frac12\log\frac{
1+k(x,y)}{1-k(x,y)}.$$
Then equality (2) assumes the form
$$g_{T(\Gamma)}(x,y)=\log k(x,y).$$
\proclaim{Corollary 3}
The function $k(x,y),$ for every fixed $y,$ is 
logarithmically
pluri\-subharmonic in $x$ on $T(\Gamma).$
\endproclaim
\proclaim{Corollary 4}
Every Teichm\"uller space $T(\Gamma)$ is  complex 
hyperconvex \RM(that means
there exists a negative plurisubharmonic function $u(x)$ 
on $T(\Gamma)$ which
tends to zero when $x$ tends to infinity\RM).
\endproclaim
The question of the complex hyperconvexity of 
Teichm\"uller spaces was stated
by M. Gromov and is solved in the affirmative for the
finite-dimensional
Teichm\"uller spaces $T(p,n)$ in \cite{9}.

Let us add here that the function $\log \kappa(x)$ 
constructed in \cite{9}
for the spaces $T(p,n)$ using the Grunsky coefficient 
inequality is related
to the Green function $g_{T(\Gamma)}(x,0)$ with strong 
inequality $\log
\kappa(x) <g_{T(\Gamma)}(x,0)$ and seems to be a minimal 
negative
plurisubharmonic function on $T(p,n)$ with prescribed 
singularity at
the origin.


\heading4. A sketch of the proof of Theorem 1\endheading

Before deriving other corollaries of Theorem 1, we shall 
give a sketch of its
proof.
For this 
we must use an alternative definition of the Green 
function proposed
by Poletskii for finite-dimensional domains (see \cite{14}).

Let $D$ be a domain in a complex Banach space $\Cal X.$
For the given points $x,y\in D,$ consider the holomorphic 
maps $f$ from the
unit disc $\Delta = \{t\in\Bbb C: |t|<1\}$ into $D$ with 
$f(0)=x$ and denote
$$\upsilon_f(x,y)=\sum_{\zeta_{j}\in 
f^{-1}(y)}k_j\log|\zeta_j|,$$
where the sum (which may be equal to $-\infty$) is taken 
over all preimages
$\zeta_j\in \Delta$ of the point $y,$ and $k_j$ is the 
multiplicity of $f$ at
$\zeta_j.$
Now put
$$\tilde g_D(x,y)=\inf\{\upsilon_f(x,y): f\in 
\text{Hol}(\Delta,D),f(0)=x\}.$$

Poletskii proved (see \cite{14}) that for the domains in 
$\bold C^n$ the
function $\tilde g_D$ coincides with the Green function 
$g_D$ defined above
by (1).
His proof uses essentially some properties of holomorphic 
functions known for
the finite-dimensional case only.
Thus we must use a different approach.

We start by proving the following
\proclaim{Lemma 1}
Let $B=B(x_0,r)$ denote the ball $\{x\in\Cal X: 
\|x-x_0\|<r\}.$
For all $(x,y)\in B\times B$ the equality $\tilde 
g_B(x,y)=g_B(x,y)$ holds.
\endproclaim

Assume now without loss of generality that $\Gamma$ acts 
discontinuously on
the upper (and lower) half plane $U=\{z:\operatorname{Im}z 
> 0\}.$
We shall apply Lemma 1 to the unit ball $M(\Gamma)$ of the 
space of Beltrami
differentials with respect to $\Gamma,$ supported in $U,$ 
with
$L^\infty$-norm.

The next considerations are different for the cases when 
$\dim
T(\Gamma)<\infty$ and $\dim T(\Gamma)=\infty,$ respectively.

In the case $\dim T(\Gamma)<\infty$ we can use the 
indicated result of
Poletskii.
It is sufficient to prove (2) for $y=0,$ the general case 
following from this
if one uses a standard right translation from $T(\Gamma)$ 
to some
$T(\Gamma').$
Applying Lemma 1 and Slodkowski's theorem on the extension 
of holomorphic
motions \cite{17}, we derive

\proclaim{Lemma 2}
For all $x\in T(\Gamma)$ the equality
$$g_{T(\Gamma)}(x,0)=\operatorname{
inf}\{g_{M(\Gamma)}(\mu,0):\Phi_{\Gamma}(\mu)=x\},$$
holds where $\Phi_\Gamma$ is the canonical projection 
$M(\Gamma)\rightarrow
T(\Gamma).$
\endproclaim

Now, using Royden's theorem about the coincidence of 
Kobayashi and 
Teichm\"uller metrics on $T(\Gamma)$ (see 
\cite{16},\cite{4}) and some
properties of the Green function, we obtain (2).

In the case $\dim T(\Gamma)=\infty$ one must approximate 
$T(\Gamma)$ by
finite-dimensional Teichm\"uller spaces following Gardiner 
\cite{5}.
\heading5. Invariant metric generated by the Green 
function\endheading

Using the negative plurisubharmonic functions, Azukawa 
\cite{1} introduced an
invariant differential (pseudo)metric on complex manifolds 
$\Cal M$ over
$\bold C^n;$ his construction was extended in \cite{13} to 
the complex Banach
manifolds.  
Let us denote this metric, determined on the tangent 
bundle $T\Cal M$ of
$\Cal M,$ by $A_{\Cal M}(x,\xi);$ here $x\in\Cal M$ and 
$\xi\in T_x\Cal M$ is
a tangent vector.

Applying Theorem 1, one can establish the following result:
\proclaim{Theorem 2}
For all spaces $T(\Gamma),$
$$A_{T(\Gamma)}(x,\xi)=F_{T(\Gamma)}(x,\xi)$$
where $F_{T(\Gamma)}(x,\xi)$ is the Finsler structure of 
$T(\Gamma).$
Consequently, the Teichm\"uller metric is nothing but the 
integrated form of
the Azukawa metric of the space $T(\Gamma).$
\endproclaim
\heading6. Plurisubharmonic functions\\ and extremal 
quasiconformal 
maps\endheading

Theorem 1 allows us to give a new characterization of 
extremal Beltrami
differentials by using plurisubharmonic functions.

Recall that a differential $\mu_0\in M(\Gamma)$ is called 
{\it extremal} if
$\|\mu_0\|=
\inf\{\|\mu\|:\Phi_{\Gamma}(\mu)=\Phi_\Gamma(\mu_0)\},$ 
i.e.,
$d_{T(\Gamma)}(0,\Phi_\Gamma(\mu_0))=d_{M(\Gamma)}(0,%
\mu_0).$
\proclaim{Theorem 3}
A Beltrami differential $\mu_0\in M(\Gamma)$ is extremal 
if and only if there
exists a logarithmically plurisubharmonic function $f\: 
T(\Gamma)\rightarrow
[0,1)$, $f(0)=0,$ such that the ratio $f(x)/\|x\|$ is 
bounded in some
neighborhood of $x=0$ and one of the following 
inequalities is satisfied{\rm:}
$$\limsup_{|t|\rightarrow 0}\frac{f\circ\Phi_\Gamma
(t\mu_{0}/\|\mu_0\|_\infty)}{|t|}=1,\qquad
f\circ\Phi_\Gamma(\mu_0)=\|\mu_0\|_\infty.$$
\endproclaim

A straightforward consequence of the above theorems is the 
following
well-known criterion for the extremality which plays an 
important role in
many applications (see, e.g.,
\cite{5, 6, 8, 15, 18}).
\proclaim{Corollary 5}
A differential $\mu_0\in M(\Gamma)$ is extremal if and 
only if
$$\sup\left\{\frac{1}{2}\left|\int\int_U\mu_0(z)\phi(z)%
\,dz\wedge d\overline
z\right|:\phi\in A_1(\Gamma),\|\phi\|=1\right\} = 
\|\mu_0\|_\infty;$$
here $A_1(\Gamma)=\{\phi\in L_1(U): \phi \ \text{is 
holomorphic in} \ U\}.$
\endproclaim
\heading7\endheading

In \cite{4} the Royden problem on the holomorphic maps 
from the disc into
the Teichm\"uller space is solved.
Theorem 1 gives another approach to this problem.

Full proofs will appear elsewhere.
\heading Acknowledgement\endheading
\par  This research began during my stay at
the University of Bielefeld, SFB 343 Diskrete Strukturen 
in der Mathematik,
and IHES, Bures-Sur-Ivette.  I am grateful for their 
support.


\enddocument